\documentclass[12pt]{article} 

\usepackage{amsmath,amsxtra,amssymb,latexsym, amscd}
\usepackage[mathscr]{eucal}

\begin{document}                        
                        \def\DATE{9/13/'03}

\title{ \huge\bf LINEAR DIFFERENCE EQUATIONS AND PERIODIC SEQUENCES 
 OVER  FINITE FIELDS 
} 

\def\1{\rule{0cm}{0cm}} \def\qd{\rule{3mm}{3mm}} \def\BB{$\bullet$}
\renewcommand{\arraystretch}{1.25}
\renewcommand{\theequation}{\thesectn.\arabic{equation}}
\def\sce{\setcounter{equation}{0}}  \newcounter{sectn} 
\newcounter{sbsect}
\def\sect#1{\addtocounter{section}{1}\sce\setcounter{sbsect}{0}%
        \renewcommand{\thesectn}{\thesection}\1\smallskip\\
        {\1\hspace{-2em}\large\bf\thesectn.\qquad #1\smallskip\par}}
\def\subsect#1{\addtocounter{sbsect}{1}\sce%
        
\renewcommand{\thesectn}{\thesection:\Alph{sbsect}}\1\smallskip\\
        {\bf\1\hspace{-1.5em}\thesectn.\qquad #1\smallskip\par}}
\newtheorem{Theorem}{THEOREM} \newtheorem{Lemma}[Theorem]{LEMMA}
\newtheorem{Corollary}[Theorem]{COROLLARY}
\def\thm#1#2{\be{Theorem}{\lb{#1} #2}} \def\LEM#1#2{\BE{Lemma}{\LB{#1} 
#2}}
\def\COR#1#2{\BE{Corollary}{\LB{#1} #2}}
\def\proof{\bigskip\noindent {\sc Proof:}\qquad}
\def\REM{\1\smallskip\par\noindent{\bf REMARK:}\qquad }
\def\qed{\hfill$\quad$\qd\medskip\\} \def\ds{\displaystyle}
\def\LB#1{\label{#1}} \def\BE#1#2{\begin{#1} #2 \end{#1}}
\def\EQ#1#2{\BE{equation}{\LB{#1} #2}} \def\ARR#1#2{\BE{array}{{#1} 
#2}}
\def\DES#1{\BE{description}{#1}} \def\QT#1{\BE{quote}{#1}}
\def\ENUM#1{\BE{enumerate}{#1}} \def\ITM#1{\BE{itemize}{#1}}
 \def\COM#1{\par\noindent{\bf COMMENT:\quad\sl #1}\par\noindent}
\def\mapsfrom{\hbox{$\;{\leftarrow}\kern-.15em{\mapstochar}\:\:$}}
\def\vv{\kern.344em{\rule[.18ex]{.075em}{1.32ex}}\kern-.344em}
\def\RE{\mbox{\rm I\kern-.21em R}} \def\CX{\mbox{\rm \vv C}}
\def\imp{\Rightarrow} \def\emb{\hookrightarrow} 
\def\wk{\rightharpoonup}
\def\rd{\dot{\1}} \def\d{\cdot} \def\+{\oplus} \def\x{\times}
\def\<{\langle} \def\>{\rangle} \def\o{\circ} \def\at#1{\Bigr|_{#1}}
\def\cd{\partial} \def\grad{\nabla} \def\L{\left} \def\R{\right}
\def\A{\mathbf  A} \def\by{\mathbf{y}} \def\bS{\mathbf{S}}
\def\I{{\cal I}} \def\A{{\mathbf A}} \def\D{{\mathbb D}}\def\bc{\mathbf{c}}
\def\X{\mathbf{X}} \def\bo{\mathbf{O}}
 \def\B{\mathbf{B}}\def\bone{\mathbf{1}}\def\bzero{\mathbf{0}}
\def\H{{\mathcal H}} \def\U{{\cal U}} \def\bc{\mathbf{c}}
\def\F{\mathbf{F}}
\def\eq{equation} \def\de{differential \eq} \def\pde{partial \de}
\def\sol{solution} \def\pb{problem} \def\bdy{boundary} 
\def\fn{function}
\def\dde{delay \de} \def\ev{eigenvalue}
\def\R{\mathbb R}
\def\Q{\mathbb Q}
\def\C{\mathbb C}
\def\N{\mathbb N}
\def\Z{\mathbb Z}
\author{
{\bf Dang Vu Giang}\\
Hanoi Institute of Mathematics\\
18 Hoang Quoc Viet, 10307 Hanoi, Vietnam\\
{\footnotesize          e-mail: $\<$dangvugiang@yahoo.com$\>$}\\
\1\\
}
\maketitle
\bigskip
\noindent {\bf Abstract.}   First, we study the linear equations in general. 
A common period is found for every solution of a linear difference equation over a finite field. 
It will help to estimate the $p$-adic modulus of polynomial roots.
Second, we focus our attention in periodic sequences over finite fields and Hamilton cycles in de Bruijn directed graph.
 


\eject
\sect{LINEAR DIFFERENCE EQUATIONS OVER P-ADIC FIELD AND CONVERGENCE }

\bigskip
\noindent
Recall  that the $p-$field ${{\mathbb{C}}_{p}}$ of  $p-$adic numbers is algebraically closed (see \cite{Robert}). Hence, every invertible $p-$adic  matrix  can be written  in the form $A=SU=US$  where $S$ is simple matrix  and $U$  is unipotent  matrix (see \cite{Goodman}).  This form is called the Jordan multiplicative decomposition of an invertible matrix $A.$  Theorems of Lie and Engel are true for Lie algebras over $p-$adic fields. Also Killing form and semisimple Lie algebra can be defined over $p-$adic fields. We are interested in the discrete dynamics of $p-$adic matrices and Haar measure of $p-$adic integers. Unlike complex matrices, $p-$adic matrices have only discrete dynamical version. (We have no sense of real and image parts of a $p-$adic number.)   Consider the difference equation ${{\mathbf{v}}_{n+1}}=A{{\mathbf{v}}_{n}}$  for  $n=0,1,2,\cdots,$ where ${{\mathbf{v}}_{0}}$ is a given vector in a $k-$dimensional vector space over ${{\mathbb{C}}_{p}}.$ Without loss of generality  we assume that $A$  is an  invertible matrix. We can write down   $A=SU=US$  where $S$ is a simple matrix  and $U$ is a unipotent  matrix (see \cite{Goodman}).  Let 
 ${{\mathbf{e}}_{1}},{{\mathbf{e}}_{2}},\cdots,{{\mathbf{e}}_{k}}$  be (linearly independent)  eigenvectors of  $S.$ Then there are eigenvalues  ${{\lambda }_{1}},{{\lambda }_{2}},\cdots ,{{\lambda }_{k}}$  (maybe identical)  of  $A$  such that  $S{{\mathbf{e}}_{1}}={{\lambda }_{1}}{{\mathbf{e}}_{1}}, \text{ }S{{\mathbf{e}}_{2}}={{\lambda }_{2}}{{\mathbf{e}}_{2}},\cdots ,S{{\mathbf{e}}_{k}}={{\lambda }_{k}}{{\mathbf{e}}_{k}}.$
Decompose  ${{\mathbf{v}}_{0}}\mathrm{=}{{\alpha }_{1}}{{\mathbf{e}}_{1}}\mathrm{+}{{\alpha }_{2}}{{\mathbf{e}}_{2}}+\cdots +{{\alpha }_{k}}{{\mathbf{e}}_{k}}$  where ${{\alpha }_{1}},{{\alpha }_{2}},\cdots ,{{\alpha }_{k}}\in {{\mathbb{C}}_{p}}.$ We have  
\[\begin{aligned}
{{\mathbf{v}}_{n}}
&={{\alpha }_{1}}{{U}^{n}}{{S}^{n}}{{\mathbf{e}}_{1}}\mathbf{+}{{\alpha }_{2}}{{U}^{n}}{{S}^{n}}{{\mathbf{e}}_{2}}+\cdots +{{\alpha }_{k}}{{U}^{n}}{{S}^{n}}{{\mathbf{e}}_{k}}\\
&={{\alpha }_{1}}\lambda _{1}^{n}{{U}^{n}}{{\mathbf{e}}_{1}}\mathbf{+}{{\alpha }_{2}}\lambda _{2}^{n}{{U}^{n}}{{\mathbf{e}}_{2}}+\cdots +{{\alpha }_{k}}\lambda _{k}^{n}{{U}^{n}}{{\mathbf{e}}_{k}},
\end{aligned}\]
which tends to zero vector  as $n\to \infty $ if  $p-$adic absolute values of  ${{\lambda }_{1}},{{\lambda }_{2}},\cdots ,{{\lambda }_{k}}$  are  less than 1.  Indeed, $U=I+N$  ($N$ is nilpotent,  $I$ is identity matrix)  so  
\[{{\lambda }^{n}}{{U}^{n}}\mathbf{e}={{\lambda }^{n}}{{\left( I+N \right)}^{n}}\mathbf{e}=\sum\limits_{\ell =0}^{k-1}{\frac{n\left( n-1 \right)\cdots \left( n-\ell +1 \right)}{\ell !}\cdot }{{\lambda }^{n}}N^\ell\mathbf{e}\]  
which is certainly tends to 0 as $n\to \infty $  if ${{\left| \lambda  \right|}_{p}}<1,$  because \[\frac{n\left( n-1 \right)\cdots \left( n-\ell +1 \right)}{\ell !}\cdot \left| \lambda  \right|_{p}^{n}\to 0\]   as  for $n\to \infty $  for  any fixed $\ell .$   We have 

\bigskip

\noindent {\bf Theorem 1.}  {\it If  the $p-$adic absolute value of every  eigenvalue of a $p-$adic matrix  $A$ is less than $1$  then every solution of ${{\mathbf{v}}_{n+1}}=A{{\mathbf{v}}_{n}}$ converges to $\mathbf 0$ as  $n\to \infty .$}

\bigskip

\noindent {\bf  Remark.}  If the $p-$adic matrix $A$ is not invertible we decompose  the $k-$ dimensional space  (over  ${{\mathbb{C}}_{p}}$)   \[V={{V}_{0}}\oplus \text{Ker}\left( A \right)\]  and set ${{A}_{0}}=A{{|}_{{{V}_{0}}}}.$  Then  ${{A}_{0}}$  is invertible and we process with ${{A}_{0}}$  and  exactly as with $A$  and  $V.$
 
\bigskip
\bigskip 

\sect{THE PERIOD OF LINEAR DIFFERENCE EQUATIONS OVER FINITE FIELDS }

\bigskip
\noindent Consider the difference equation ${{\mathbf{v}}_{n+1}}=A{{\mathbf{v}}_{n}}$  for  $n=0,1,2,\cdots ,$ where ${{\mathbf{v}}_{0}}$ is a given vector in a $k-$dimensional vector space over ${{\mathbb{F}}_{q}}$ ($q$ is a power of a prime $p$). First we assume that the matrix $A$  is invertible. Then  $A=SU=US$  where $S$  is simple matrix  and $U$  is unipotent  matrix over an extension of the field ${{\mathbb{F}}_{q}}$ (see Goodman et al).  We have  
\[\begin{aligned}
{{\mathbf{v}}_{n}}
&={{\alpha }_{1}}{{U}^{n}}{{S}^{n}}{{\mathbf{e}}_{1}}\mathrm{+}{{\alpha }_{2}}{{U}^{n}}{{S}^{n}}{{\mathbf{e}}_{2}}+\cdots +{{\alpha }_{k}}{{U}^{n}}{{S}^{n}}{{\mathbf{e}}_{k}}\\
&={{\alpha }_{1}}\lambda _{1}^{n}{{U}^{n}}{{\mathbf{e}}_{1}}\mathrm{+}{{\alpha }_{2}}\lambda _{2}^{n}{{U}^{n}}{{\mathbf{e}}_{2}}+\cdots +{{\alpha }_{k}}\lambda _{k}^{n}{{U}^{n}}{{\mathbf{e}}_{k}},
\end{aligned}\]
where ${{\mathbf{e}}_{1}},{{\mathbf{e}}_{2}},\cdots ,{{\mathbf{e}}_{k}}$ are  (linearly independent)  eigenvectors of $S$ and $\lambda _{1}^{{}},\lambda _{2}^{{}},$  $\cdots ,\lambda _{k}^{{}}$  are eigenvalues of $A$(in a finite extension of ${{\mathbb{F}}_{q}}$). Moreover, $U=N+I$ ($N$ is nilpotent, $I$  is identity matrix) so  
  \[{{\lambda }^{n}}{{U}^{n}}\mathbf{e}={{\lambda }^{n}}{{\left( I+N \right)}^{n}}\mathbf{e}=\sum\limits_{\ell =0}^{k-1}{\frac{n\left( n-1 \right)\cdots \left( n-\ell +1 \right)}{\ell !}\cdot }{{\lambda }^{n}}N^\ell\mathbf{e}.\]
A lot is known about the periods of vector sequences over ${\mathbb F}_q$ defined by a recurrence ${\mathbf v}_{n+1}=A{\mathbf v}_n$ with a matrix $A$, see for example \cite{Niederreiter2} or Chapter 10.1 of \cite{Niederreiter1}.
However, the explicit formula of the period is found for irreducible characteristic polynomials  only. 
Here we give the explicit formula of the period of  ${\mathbf v}_n$ not only in that case but for arbitrary characteristic polynomial. Our method is based on  Lucas' congruence and writing the dimension $k$ in $p-$adic base. Indeed, 
the period of  $\left( \begin{matrix}
   n  \\
   \ell  \\
\end{matrix} \right)$ mod $p$ is $p$ if $\ell<p$. Clearly, ${{\lambda }_{j}}$  is root of the characteristic equation det$\left( A-\lambda I \right)=0$ (in an extension of ${{\mathbb{F}}_{q}}$). Consider the characteristic equation det$\left( A-\lambda I \right)=0$ in an extension of ${{\mathbb{F}}_{q}}$ and let $\left\{ {{\varphi }_{j}} \right\}_{j=1}^{N}\subseteq {{\mathbb{F}}_{q}}\left[ x \right]$ are irreducible factors (of the characteristic polynomial of the matrix $A$).  If  ${{k}_{j}}$  is degree of  ${{\varphi }_{j}}$  then ${{\varphi }_{j}}$  is  a divisor of ${{x}^{{{q}^{{{k}_{j}}}}}}-x$  (see \cite{Lidl} for more details). Let $Q$  be the smallest multiplier of  ${{q}^{{{k}_{j}}}}-1$ then every root $\lambda $  of the characteristic equation det$\left( A-\lambda I \right)=0$ satisfies ${{\lambda }^{Q}}=1$  so the period of  ${{\mathbf{v}}_{n+1}}=A{{\mathbf{v}}_{n}}$ is $Q$  if there is no duplicated root ($U=I$ in this case).  If there is a duplicated root, we should  make an intention about the periodicity of  sequences $\left( \begin{matrix}
   n  \\
   \ell  \\
\end{matrix} \right)$ of $n$ in  ${{\mathbb{F}}_{q}}.$  By Lucas' congruence 
\[{n \choose j}\equiv {n_0 \choose j_0} \cdots {n_r \choose j_r}\mod p\]
 if  $n=n_0+\cdots+n_rp^{r}+\cdots$ and $j=j_0+\cdots+j_rp^{r}$ with $0\le n_i,j_i<p$. Hence, if we write $k-1={{\kappa }_{0}}+\cdots+{{\kappa }_{r}}{{p}^{r}}$ then the period of sequences $\left( \begin{matrix}
   n  \\
   \ell  \\
\end{matrix} \right)$  $\left( 0\le \ell \le k-1 \right)$ is ${{p}^{r+1}}.$

\bigskip

\par\noindent{\bf  Theorem 2. } {\it   Assume that $A$ is $k\times k$  invertible matrix over ${{\mathbb{F}}_{q}}.$  Let ${{k}_{1}}<{{k}_{2}}<\cdots <{{k}_{N}}$  be degrees of irreducible factors of the characteristic polynomial of $A$ and let $Q=Q\left( A \right)$ be the smallest multiplier of  ${{q}^{{{k}_{j}}}}-1.$ Then every solution of the linear system ${{\mathbf{v}}_{n+1}}=A{{\mathbf{v}}_{n}}$  over  ${{\mathbb{F}}_{q}}$ is $Q-$periodic if the multiplicity of every irreducible factor is $1$. Specially if the characteristic polynomial of $A$ is irreducible then $q^k-1$ is a common period of every 
solution of the linear system ${{\mathbf{v}}_{n+1}}=A{{\mathbf{v}}_{n}}$  over  ${{\mathbb{F}}_{q}}$.
If  the multiplicity of a factor is larger than $1$  then every solution of 
${{\mathbf{v}}_{n+1}}=A{{\mathbf{v}}_{n}}$  is $Qp^{r+1}-$periodic, where $p$ is the characteristics  of  ${{\mathbb{F}}_{q}}$  and the exponent $r$ is determined by writing 
$k-1={{\kappa }_{0}}+\cdots+{{\kappa }_{r}}{{p}^{r}}$ with $0\le {{\kappa }_{i}}<p$ and $0<{{\kappa }_{r}}<p.$ 
}
 
\bigskip

\par\noindent{\bf  Remark. }  If the matrix  $A$ is not invertible we decompose  the $k-$ dimensional space (over  ${{\mathbb{F}}_{q}}$) $V={{V}_{0}}\oplus \text{Ker}\left( A \right)$  and set ${{A}_{0}}=A{{|}_{{{V}_{0}}}}.$  Then ${{A}_{0}}$  is invertible and we process with ${{A}_{0}}$    exactly as with $A.$ The system ${{\mathbf{v}}_{n+1}}=A{{\mathbf{v}}_{n}}$  is $\rho -$like periodic, that is there is a period $T\in \left\{ Q\left( {{A}_{0}} \right),pQ\left( {{A}_{0}} \right) \right\}$  such that ${{\mathbf{v}}_{n+T}}={{\mathbf{v}}_{n}}$  for all but finite $n$. Specially, if $A$  is nilpotent,  ${{\mathbf{v}}_{n}}=\mathbf{0}$ for $n=k,k+1,\cdots $  and the period  $T=1.$
  
\bigskip

\par\noindent{\bf   
Example 1. } Let ${{a}_{n}}=n$ for $n=0,1,\cdots ,p-1$ and ${{a}_{n}}={{a}_{n-1}}+{{a}_{n-p}}$   for   $n=p,p+1,\cdots .$  In the finite field ${{\mathbb{F}}_{p}}$ we have ${{a}_{n}}={{a}_{n+{{p}^{2}}-1}}.$ To see this, solve the equation ${{\lambda }^{p}}={{\lambda }^{p-1}}+1$  in an extension of  ${{\mathbb{F}}_{p}}.$ We have ${{\left( \lambda -1 \right)}^{p}}={{\lambda }^{p-1}}$  and consequently, ${{\left( \lambda -1 \right)}^{p+1}}={{\lambda }^{p-1}}\left( \lambda -1 \right)$ $=1.$ Therefore, ${{\lambda }^{{{p}^{2}}-1}}={{\left( \lambda -1 \right)}^{p\left( p+1 \right)}}=1.$ Hence ${{a}_{n}}={{a}_{n-1}}+{{a}_{n-p}}$ is periodic in ${{\mathbb{F}}_{p}}$ with period ${{p}^{2}}-1.$ This sequence is not convergent in the $p-$adic field, since ${{\left| {{a}_{n+\ell \left( {{p}^{2}}-1 \right)}} \right|}_{p}}=1$  for  all $\ell \in \mathbb{N}$  and  $n=1,\cdots ,p-1$ so ${{\left| {{a}_{n+\ell \left( {{p}^{2}}-1 \right)+p}}-{{a}_{n+\ell \left( {{p}^{2}}-1 \right)+p-1}} \right|}_{p}}={{\left| {{a}_{n+\ell \left( {{p}^{2}}-1 \right)}} \right|}_{p}}=1$  which does not converge to 0.  Hence,  ${{a}_{n}}={{a}_{n-1}}+{{a}_{n-p}}$ is not convergent in the $p-$adic field and at least one root of  ${{\lambda }^{p}}={{\lambda }^{p-1}}+1$ has $p-$adic absolute value $\ge 1.$
 
 \bigskip

\par\noindent{\bf   
Example 2. } 
Consider the equation ${{x}_{n}}={{x}_{n-4}}+{{x}_{n-5}}$ over ${{\mathbb{F}}_{2}}.$ The characteristic polynomial is ${{\lambda }^{5}}+\lambda +1=\left( {{\lambda }^{2}}+\lambda +1 \right)\left( {{\lambda }^{3}}+\lambda^2 +1 \right).$ The smallest multiplier of
${{2}^{2}}-1=3$  and  ${{2}^{3}}-1=7$ is 21. Hence a period of this equation is 21. This example is given by \cite{Berg}. 

 \bigskip

\par\noindent{\bf   
Example 3. } 
 Consider ${{x}_{n}}={{x}_{n-3}}+{{x}_{n-5}}$  over  ${{\mathbb{F}}_{2}}.$ The characteristic polynomial is ${{\lambda }^{5}}+{{\lambda }^{2}}+1$   irreducible over  ${{\mathbb{F}}_{2}}.$ Hence the minimal period of  any non-zero solution of this equation is ${{2}^{5}}-1=31.$ This example is given by \cite{Berg}. 
  \bigskip

\par\noindent{\bf   
Example 4. } 
Consider the equation ${{x}_{n}}+{{x}_{n+1}}+{{x}_{n+2}}+{{x}_{n+4}}=0$  over ${{\mathbb{F}}_{2}}.$  The characteristic polynomial is  ${{\lambda }^{4}}+{{\lambda }^{2}}+\lambda +1=\left( \lambda +1 \right)\left( {{\lambda }^{3}}+{{\lambda }^{2}}+1 \right).$ The smallest  multiplier of  ${{2}^{1}}-1=1$  and  ${{2}^{3}}-1=7$ is 7. Hence the minimal period of any non-zero solution of this equation is 7. This example is given by \cite{Berg}. 
 
 \bigskip

\par\noindent{\bf   
Example 5. } 
Consider  the  Fibonacci  sequence ${{x}_{n}}={{x}_{n-1}}+{{x}_{n-2}}$  over ${{\mathbb{F}}_{q}}.$ The characteristic polynomial is ${{\lambda }^{2}}-\lambda -1.$  If this polynomial is irreducible over ${{\mathbb{F}}_{q}}$ then the Fibonacci sequence is $\left( {{q}^{2}}-1 \right)-$periodic. Otherwise, it is $\left( q-1 \right)-$ or $\left( q-1 \right)p-$periodic. For example, ${{\lambda }^{2}}-\lambda -1$  is irreducible  over ${{\mathbb{F}}_{2}}$ and ${{\mathbb{F}}_{3}}$ but ${{\lambda }^{2}}-\lambda -1={{\left( \lambda -3 \right)}^{2}}$  over  ${{\mathbb{F}}_{5}}.$ Hence, the Fibonacci sequence is $3-$periodic over ${{\mathbb{F}}_{2}}$,   $8-$periodic over ${{\mathbb{F}}_{3}}$  and 
 $20-$periodic over  ${{\mathbb{F}}_{5}}.$ In ${{\mathbb{F}}_{4}},$ ${{\lambda }^{2}}-\lambda -1$  has two distinct roots so the Fibonacci sequence  is $3-$periodic over  ${{\mathbb{F}}_{4}}.$ This example is solved in part by \cite{Ladas}. 

\bigskip
\bigskip 

\sect{MINIMAL POLYNOMIAL AND MINIMAL PERIOD OF PERIODIC SEQUENCES}

\bigskip

\par\noindent  Let $S=\left\{ {{s}_{0}},{{s}_{1}},\cdots  \right\}$ be an arbitrary sequence over ${{\mathbb{F}}_{q}}$. Let $\mathbb{E}$ denote the shift operator on the linear space of sequences. This means that 
 $\mathbb{E}\left( S \right)$ is  the sequence whose $i$th term is ${{s}_{i+1}}.$ Moreover, ${{\mathbb{E}}^{k}}\left( S \right)$ is denoted the sequence whose $i$th term is ${{s}_{i+k}}.$ Two sequences $S=\left\{ {{s}_{0}},{{s}_{1}},\cdots  \right\}$ and $T=\left\{ {{t}_{0}},{{t}_{1}},\cdots  \right\}$ are equivalent if ${{\mathbb{E}}^{k}}\left( S \right)=T$ for some integer $k$.
Suppose that 
	\[{{s}_{i+n}}+{{c}_{n-1}}{{s}_{i+n-1}}+\cdots +{{c}_{0}}{{s}_{i}}=0\]
for some elements ${{c}_{0}},{{c}_{1}},\cdots ,{{c}_{n-1}}\in {{\mathbb{F}}_{q}}$ and for all $i=0,1,\cdots .$ Then the polynomial ${{x}^{n}}+{{c}_{n-1}}{{x}^{n-1}}+\cdots +{{c}_{0}}$ is called a characteristic polynomial of $S=\left\{ {{s}_{0}},{{s}_{1}},\cdots  \right\}$. For a periodic sequence $S=\left\{ {{s}_{0}},{{s}_{1}},\cdots  \right\}$ with a period $N$ we have ${{s}_{i+N}}-{{s}_{i}}=0$, so ${{x}^{N}}-1$ is a characteristic polynomial of $S=\left\{ {{s}_{0}},{{s}_{1}},\cdots  \right\}$. A characteristic polynomial of minimal degree is called the {\bf minimal polynomial } of $S=\left\{ {{s}_{0}},{{s}_{1}},\cdots  \right\}$.  Any periodic non-zero  sequence has one and only one minimal polynomial. The degree of  the minimal polynomial is called the {\bf complexity}  of the periodic sequence and denoted by $c\left( S \right).$ Clearly, $c\left( S \right)\le $ minimal period of $S=\left\{ {{s}_{0}},{{s}_{1}},\cdots  \right\}.$
Moreover, the minimal polynomial divides any characteristic polynomial of $S=\left\{ {{s}_{0}},{{s}_{1}},\cdots  \right\}$. Hence, if a characteristic polynomial of a periodic sequence  is irreducible then it is itself  the minimal polynomial of that sequence.
If $N$ is a period of $S=\left\{ {{s}_{0}},{{s}_{1}},\cdots  \right\}$ then the minimal polynomial of $S=\left\{ {{s}_{0}},{{s}_{1}},\cdots  \right\}$ is 
\[\frac{{{x}^{N}}-1}{\text{gcd}\left( {{x}^{N}}-1,{{s}_{0}}+{{s}_{1}}x+\cdots +{{s}_{N-1}}{{x}^{N-1}} \right)}.\]
Therefore, 
\[c\left( S \right)=N-\deg \text{gcd}\left( {{x}^{N}}-1,{{s}_{0}}+{{s}_{1}}x+\cdots +{{s}_{N-1}}{{x}^{N-1}} \right).\]
For example, the linear complexity of the Fibonacci sequence $\{1,1,2,3,\cdots\}$ over $\mathbb F_5$ 
is 2 and its minimal period is 20.

\bigskip
\par\noindent{\bf   Theorem 3.} {\it 
The minimal period of a non-zero sequence $S=\left\{ {{s}_{0}},{{s}_{1}},\cdots  \right\}$ is a power of $p$ iff the minimal polynomial of $S=\left\{ {{s}_{0}},{{s}_{1}},\cdots  \right\}$ is a power of $\left( x-1 \right)$. Therefore, if a characteristic polynomial of $S=\left\{ {{s}_{0}},{{s}_{1}},\cdots  \right\}$ is of the form ${{\left( x-1 \right)}^{k}}$ then the minimal period of $S=\left\{ {{s}_{0}},{{s}_{1}},\cdots  \right\}$ is ${{p}^{r+1}}$ with ${{p}^{r}}< c\left( S \right)\le {{p}^{r+1}}$, unless $c\left( S \right)=1$, in which case  
$S=\left\{ {{s}},{{s}},\cdots  \right\}$ and the minimal period of $S$ is $1$. 
Recall that $p$ is the characteristics of the finite field ${{\mathbb{F}}_{q}}$.
}
\bigskip
\par\noindent{\sl
Proof: }
If ${{p}^{r+1}}$ is the minimal period of $S=\left\{ {{s}_{0}},{{s}_{1}},\cdots  \right\}$ then 
\[ {{x}^{{{p}^{r+1}}}}-1 ={{\left( x-1 \right)}^{{{p}^{r+1}}}}\]
 is a characteristic polynomial of $S=\left\{ {{s}_{0}},{{s}_{1}},\cdots  \right\}$. Therefore, ${{\left( x-1 \right)}^{c\left( S \right)}}$ is the minimal polynomial of $S=\left\{ {{s}_{0}},{{s}_{1}},\cdots  \right\}$.
Conversely, if the minimal polynomial of $S=\left\{ {{s}_{0}},{{s}_{1}},\cdots  \right\}$ is ${{\left( x-1 \right)}^{c\left( S \right)}}$
 then the minimal period of $S=\left\{ {{s}_{0}},{{s}_{1}},\cdots  \right\}$ is a power of $p$.
Indeed, if  $c\left( S \right)=1$ then $S=\left\{ {{s}},{{s}},\cdots  \right\}$ and the minimal period is 1. For   $c\left( S \right)>1$ choose an integer $r$ such that ${{p}^{r}}< c\left( S \right)\le {{p}^{r+1}}$ (this integer $r$ is unique). Then ${{\left( x-1 \right)}^{{{p}^{r+1}}}}={{x}^{{{p}^{r+1}}}}-1$  is a characteristic polynomial of $S=\left\{ {{s}_{0}},{{s}_{1}},\cdots  \right\}$ so ${{p}^{r+1}}$ is a period of $S=\left\{ {{s}_{0}},{{s}_{1}},\cdots  \right\}.$ Moreover, ${{p}^{r+1}}$ is the minimal period of $S=\left\{ {{s}_{0}},{{s}_{1}},\cdots  \right\}$ because  $c\left( S \right)>{{p}^{r}}.$  
The proof is complete.

\bigskip
\noindent {\bf  Remark.}
These results can be found in \cite{Paterson} with more details. Recall that the weight of a periodic sequence $S=\left\{ {{s}_{0}},{{s}_{1}},\cdots  \right\}$ with minimal period $t$ is the sum \[{{s}_{1}}+{{s}_{2}}+\cdots +{{s}_{t}}.\]

\par\noindent{\bf   Proposition.} {\it  If the minimal period of $S=\left\{ {{s}_{0}},{{s}_{1}},\cdots  \right\}$ is ${{p}^{r+1}}$ then $c\left( S \right)={{p}^{r+1}}$ if and only if  the weight of $S=\left\{ {{s}_{0}},{{s}_{1}},\cdots  \right\}$  is non-zero. }

\bigskip
\par\noindent{\sl
Proof: }  If \[{{s}_{1}}+{{s}_{2}}+\cdots +{{s}_{{{p}^{r+1}}}}=0\] then \[{{s}_{n+1}}+{{s}_{n+2}}+\cdots +{{s}_{n+{{p}^{r+1}}}}=0\] for every $n=0,1,\cdots, $ so ${{x}^{{{p}^{r+1}}-1}}+{{x}^{{{p}^{r+1}}-2}}+\cdots +1$ is a characteristic polynomial and consequently,  $c\left( S \right)<{{p}^{r+1}}.$  If 
	\[{{s}_{1}}+{{s}_{2}}+\cdots +{{s}_{{{p}^{r+1}}}}\ne 0\]
 then 
	\[{{s}_{n+1}}+{{s}_{n+2}}+\cdots +{{s}_{n+{{p}^{r+1}}}}\ne 0\]
 for every $n=0,1,\cdots,$ so ${{x}^{{{p}^{r+1}}-1}}+{{x}^{{{p}^{r+1}}-2}}+\cdots +1={{\left( x-1 \right)}^{{{p}^{r+1}}-1}}$ is not a characteristic polynomial and consequently,  $c\left( S \right)={{p}^{r+1}}.$ 
The proof is complete.

\bigskip
\noindent {\bf  Remark.}  This proposition is proved in \cite {Paterson2} by other way.

\bigskip
\par\noindent{\bf   Proposition.} {\it Let ${{m}_{1}},{{m}_{2}},\cdots ,{{m}_{h}}$ be minimal polynomials of  periodic sequences ${{S}_{1}},{{S}_{2}},\cdots ,{{S}_{h}}$ and assume that they are pairwise relatively prime.
 Then the minimal polynomial of ${{S}_{1}}+{{S}_{2}}+\cdots +{{S}_{h}}$ is ${{m}_{1}}{{m}_{2}}\cdots {{m}_{h}}.$
}
\bigskip
\par\noindent{\sl
Proof: } By complete induction, it is enough to consider the case $h=2.$ Let $N$ be a common period of ${{S}_{1}}$ and ${{S}_{2}}.$ Then 
\[{{m}_{1}}\left( x \right)=\frac{{{x}^{N}}-1}{\text{gcd}\left( {{x}^{N}}-1,{{S}_{1}}\left( x \right) \right)}\text{  and  }{{m}_{2}}\left( x \right)=\frac{{{x}^{N}}-1}{\text{gcd}\left( {{x}^{N}}-1,{{S}_{2}}\left( x \right) \right)}.\]
On the other hand, we can write  ${{x}^{N}}-1={{f}_{1}}{{m}_{1}}={{f}_{2}}{{m}_{2}}$ and ${{S}_{1}}={{f}_{1}}{{g}_{1}},$  ${{S}_{2}}={{f}_{2}}{{g}_{2}}$  with $\text{gcd}\left\{ {{m}_{1}},{{g}_{1}} \right\}=\text{gcd}\left\{ {{m}_{2}},{{g}_{2}} \right\}=1.$ But ${{m}_{1}}$and ${{m}_{2}}$ are co-prime so ${{f}_{1}}={{m}_{2}}\delta $ and ${{f}_{2}}={{m}_{1}}\delta $ and we have ${{x}^{N}}-1={{m}_{1}}{{m}_{2}}\delta $ and ${{S}_{1}}+{{S}_{2}}=\left( {{m}_{2}}{{g}_{1}}+{{m}_{1}}{{g}_{2}} \right)\delta .$ On the other hand the minimal polynomial of ${{S}_{1}}+{{S}_{2}}$ is 
\[\begin{aligned}
 \frac{{{x}^{N}}-1}{\text{gcd}\left( {{x}^{N}}-1,{{S}_{1}}\left( x \right)+{{S}_{2}}\left( x \right) \right)}
  &=\frac{{{m}_{1}}{{m}_{2}}\delta }{\text{gcd}\left( {{m}_{1}}{{m}_{2}}\delta ,\left( {{m}_{2}}{{g}_{1}}+{{m}_{1}}{{g}_{2}} \right)\delta  \right)} \\ 
 & ={{m}_{1}}{{m}_{2}}.  
\end{aligned}\]
The proof is now complete.

\bigskip
\par\noindent{\bf 	
Definition.}  For a monic polynomial $f\in {{\mathbb{F}}_{q}}\left[ x \right]$ we define $\mathcal M\left( f \right)$ to be the vector space of periodic sequences having $f$ as a characteristic polynomial. Then the dimension of $\mathcal M\left( f \right)$ over $ {{\mathbb{F}}_{q}}$  is deg$\left( f \right).$	Moreover, $\mathcal M\left( f \right)$ is invariant under the shift operator $\mathbb{E}$.
	
\bigskip
\par\noindent{\bf Proposition.}  {\it Let ${{f}_{1}},{{f}_{2}},\cdots ,{{f}_{h}}\in {{\mathbb{F}}_{q}}\left[ x \right]$ be monic polynomials and pairwise relatively prime. Then $\mathcal M\left( {{f}_{1}}{{f}_{2}}\cdots {{f}_{h}} \right)=\mathcal  M\left( {{f}_{1}} \right)\oplus \mathcal M\left( {{f}_{2}} \right)\oplus \cdots \oplus \mathcal M\left( {{f}_{h}} \right).$
}

\bigskip
\par\noindent{\sl
Proof: } By complete induction, it is enough to consider the case $h=2.$ If $S\in M\left( {{f}_{1}}{{f}_{2}} \right)$ then ${{f}_{2}}\left( \mathbb{E} \right)S\in\mathcal M\left( {{f}_{1}} \right)$ and ${{f}_{1}}\left( \mathbb{E} \right)S\in\mathcal M\left( {{f}_{2}} \right).$ Recall that $\mathbb{E}$ is the shift operator on the vector spaces of periodic sequences. On the other hand, since ${{f}_{1}},{{f}_{2}}$ are relatively prime, there are
${{g}_{1}},{{g}_{2}}\in {{\mathbb{F}}_{q}}\left[ x \right]$ such that
  ${{f}_{1}}{{g}_{1}}+{{f}_{2}}{{g}_{2}}=1$. Therefore, $S={{g}_{1}}{{f}_{1}}\left( \mathbb{E} \right)S+{{g}_{2}}{{f}_{2}}\left( \mathbb{E} \right)S\in\mathcal M\left( {{f}_{1}} \right)\oplus\mathcal M\left( {{f}_{2}} \right).$
The proof is now complete.

\bigskip
\par\noindent{\bf Proposition.}  {\it Let $m={{m}_{1}}{{m}_{2}}\cdots {{m}_{h}}$ be the minimal polynomial of a periodic sequence  $S$, where ${{m}_{1}},{{m}_{2}},\cdots ,{{m}_{h}}$ are monic polynomials and  pairwise relatively prime. Then there are unique periodic sequences  ${{S}_{1}},{{S}_{2}},\cdots ,{{S}_{h}}$ with minimal polynomials ${{m}_{1}},{{m}_{2}},\cdots ,{{m}_{h}}$, respectively, such that \[{{S}_{1}}+{{S}_{2}}+\cdots +{{S}_{h}}=S.\]
}

\par\noindent{\sl
Proof: } By the above proposition, $S={{S}_{1}}+{{S}_{2}}+\cdots +{{S}_{h}}$ with ${{S}_{j}}\in\mathcal  M\left( {{m}_{j}} \right)$ for $j=1,2,\cdots ,h.$ Let ${{\mu }_{1}},{{\mu }_{2}},\cdots ,{{\mu }_{h}}$ be the minimal polynomials of ${{S}_{1}},{{S}_{2}},\cdots ,{{S}_{h}},$ respectively. Then  ${{\mu }_{j}}|{{m}_{j}}$ for $j=1,2,\cdots ,h$ and the minimal polynomial of $S={{S}_{1}}+{{S}_{2}}+\cdots +{{S}_{h}}$ is $m={{m}_{1}}{{m}_{2}}\cdots {{m}_{h}}={{\mu }_{1}}{{\mu }_{2}}\cdots {{\mu }_{h}}.$ Therefore, ${{\mu }_{j}}={{m}_{j}}$ for $j=1,2,\cdots ,h$ and the minimal polynomials of ${{S}_{1}},{{S}_{2}},\cdots ,{{S}_{h}}$ are ${{m}_{1}},{{m}_{2}},\cdots ,{{m}_{h}},$ respectively. The proof is now complete.

\bigskip
\par\noindent{\bf 	
Definition.}  Let  $\sigma :{{\mathbb{F}}_{q}}\left[ x \right]\to {{\mathbb{F}}_{q}}\left[ x \right]$ be the linear mapping
defined by $\sigma \left( \sum{{{c}_{j}}{{x}^{j}}} \right)=\sum{c_{j}^{p}{{x}^{j}}}$. Recall that  $p$ is the characteristics of the finite field ${{\mathbb{F}}_{q}}$. 
Then $\sigma$ is also multiplicative map that is 
$\sigma \left( fg \right)=\sigma \left( f \right)\sigma \left( g \right).$
Moreover, note that the set
$\left\{ a\in {{\mathbb{F}}_{q}}:{{a}^{p}}=a \right\}$ is a field of $p$ elements so $\sigma \left( f \right)=f$ iff $f\in {{\mathbb{F}}_{p}}\left[ x \right]$ i.e., the coefficients of $f$ are in the prime field ${{\mathbb{F}}_{p}}.$
Let $k\left( f \right)$ denote  the smallest integer $\ell $ such that ${{\sigma }^{\ell }}\left( f \right)=f.$  Clearly, the set
$\left\{ a\in {{\mathbb{F}}_{q}}:{{a}^{{{p}^{k\left( f \right)}}}}=a \right\}$ is a field of ${{p}^{k\left( f \right)}}$ elements so $f\in {{\mathbb{F}}_{{{p}^{k\left( f \right)}}}}\left[ x \right].$ Moreover,  if $q={{p}^{n}}$ then $k\left( f \right)|n$.

\bigskip
\par\noindent{\bf   Proposition.} {\it
A polynomial $f\in {{\mathbb{F}}_{q}}\left[ x \right]$ is irreducible if and only if  $\sigma \left( f \right)$ is irreducible.
}
\bigskip
\par\noindent{\sl
Proof: }   If $f=uv$ is reducible then $\sigma \left( f \right)=\sigma \left( u \right)\sigma \left( v \right)$ is also reducible. Conversely, if $\sigma \left( f \right)=\mu \nu $ is reducible then $f={{\sigma }^{\ell }}\left( f \right)={{\sigma }^{\ell -1}}\left( \mu \nu  \right)={{\sigma }^{\ell -1}}\left( \mu  \right){{\sigma }^{\ell -1}}\left( \nu  \right)$ is also reducible. The proof is now complete.

\bigskip
\par\noindent{\bf   Proposition.} {\it If $f\left( x \right)\in {{\mathbb{F}}_{q}}\left[ x \right]$ is a monic irreducible polynomial then  
\[g=f\sigma \left( f \right){{\sigma }^{2}}\left( f \right)\cdots {{\sigma }^{k\left( f \right)-1}}\left( f \right)\]
 is in ${{\mathbb{F}}_{p}}\left[ x \right]$
 and irreducible over the prime field ${{\mathbb{F}}_{p}}.$ Recall that $k\left( f \right)$ denotes the smallest integer $\ell $ such that ${{\sigma }^{\ell }}\left( f \right)=f.$
}

\bigskip
\par\noindent{\sl
Proof: }   Note that $\sigma \left( g \right)=\sigma \left( f \right){{\sigma }^{2}}\left( f \right)\cdots {{\sigma }^{k\left( f \right)}}\left( f \right)=g$ so $g\in {{\mathbb{F}}_{p}}\left[ x \right]$. Moreover, 
we can decompose $f\left( x \right)=\left( x-{{\alpha }_{1}} \right)\left( x-{{\alpha }_{2}} \right)\cdots \left( x-{{\alpha }_{n}} \right)$ where ${{\alpha }_{j}}\in {{\mathbb{F}}_{{{q}^{n}}}}$ for  $j=1,2,\cdots ,n.$
Then $f$ is the minimal polynomial of ${{\alpha }_{1}}$. 
Let $\gamma $ denote the minimal polynomial of ${{\alpha }_{1}}$ over the prime field ${{\mathbb{F}}_{p}}.$ Then $\gamma |g$ and $\gamma $ is  irreducible over the prime field ${{\mathbb{F}}_{p}}.$ Moreover, $f|\gamma$
so ${{\alpha }_{j}},\alpha _{j}^{p},\alpha _{j}^{{{p}^{2}}},\cdots $ are roots of $\gamma $ for  $j=1,2,\cdots ,n$.
Therefore, $\gamma =g$ and consequently, $g$ irreducible over the prime field ${{\mathbb{F}}_{p}}.$
 The proof is now complete.

\bigskip
\par\noindent{\bf   Proposition.} {\it If  the minimal polynomial $m$ of a  periodic sequence $S=\left\{ {{s}_{0}},{{s}_{1}},\cdots  \right\}$ is irreducible then
the minimal polynomial of  ${{S}^{p}}=\left\{ s_{0}^{p},s_{1}^{p},\cdots  \right\}$ is $\sigma \left( m \right)$, which is also  irreducible.
} 

\bigskip
\par\noindent{\sl
Proof: } Since  $m$ is   irreducible, $\sigma \left( m \right)$  is irreducible.
On the other hand if 
	\[{{s}_{i+n}}+{{c}_{n-1}}{{s}_{i+n-1}}+\cdots +{{c}_{0}}{{s}_{i}}=0\] 
then 
	\[{{s}_{i+n}^p}+{{c}_{n-1}^p}{{s}^p_{i+n-1}}+\cdots +{{c}^p_{0}}{{s}^p_{i}}=0.\]
Therefore, $\sigma \left( m \right)$  is a characteristic polynomial of  ${{S}^{p}}=\left\{ s_{0}^{p},s_{1}^{p},\cdots  \right\}$. But $\sigma \left( m \right)$  is   irreducible so it should be the minimal polynomial of  ${{S}^{p}}=\left\{ s_{0}^{p},s_{1}^{p},\cdots  \right\}$.
The proof is now complete.

\bigskip
\par\noindent{\bf  Theorem 4.} {\it
 If the minimal polynomial $f$ of $S=\left\{ {{s}_{0}},{{s}_{1}},\cdots  \right\}$ is irreducible then the minimal polynomial  of  $S+{{S}^{p}}+{{S}^{{{p}^{2}}}}+\cdots +{{S}^{{{p}^{k\left( f \right)-1}}}}$ is $f\sigma \left( f \right){{\sigma }^{2}}\left( f \right)\cdots {{\sigma }^{k\left( f \right)-1}}\left( f \right)$ which is in ${{\mathbb{F}}_{p}}\left[ x \right]$ and irreducible over the prime field ${{\mathbb{F}}_{p}}.$ Hence, the comlexity of the sum $S+{{S}^{p}}+{{S}^{{{p}^{2}}}}+\cdots +{{S}^{{{p}^{k\left( f \right)-1}}}}$ is
 } $c(S)k(f)$.

\bigskip
\par\noindent{\sl
Proof: } Since $f$ is irreducible, $\sigma \left( f \right),{{\sigma }^{2}}\left( f \right),\cdots ,{{\sigma }^{k\left( f \right)-1}}\left( f \right)$ are all irreducible. But they are the minimal polynomials  of  $S,{{S}^{p}},{{S}^{{{p}^{2}}}},\cdots ,{{S}^{{{p}^{k\left( f \right)-1}}}},$ respectively. Hence, $f,\sigma \left( f \right),{{\sigma }^{2}}\left( f \right),\cdots ,{{\sigma }^{k\left( f \right)-1}}\left( f \right)$ are pairwise relatively prime and the minimal polynomial  of the sum $S+{{S}^{p}}+{{S}^{{{p}^{2}}}}+\cdots +{{S}^{{{p}^{k\left( f \right)-1}}}}$ is $f\sigma \left( f \right){{\sigma }^{2}}\left( f \right)\cdots {{\sigma }^{k\left( f \right)-1}}\left( f \right)$. 
 The proof is now complete.
 
  \bigskip
\par\noindent{\bf   Proposition.} {\it If  the minimal polynomial $f$ of a periodic sequence over 
${{\mathbb{F}}_{q}}$
is irreducible then} ${{p}^{k\left( f \right)\text{deg}\left( f \right)}}-1$  {\it is a period of this sequence.
}

\bigskip
\par\noindent{\sl
Proof: } Apply Theorem 2 in the special case where the characteristic polynomial of the matrix $A$ is irreducible and replace $q$ by ${{p}^{k\left( f \right)}}$. This application is possible because  $f\in {{\mathbb{F}}_{{{p}^{k\left( f \right)}}}}\left[ x \right]$ and $f$ is irreducible over there.

\bigskip
\par\noindent{\bf 	
Definition. The order of a polynomial}  $f\in {{\mathbb{F}}_{q}}\left[ x \right]$ with $f\left( 0 \right)\ne 0$ is the smallest positive integer $\ell $ such that $f$ divides ${{x}^{\ell }}-1$. Let ord$\left( f \right)$ denote the order of $f.$ If $f$ is the minimal polynomial of a periodic sequence over ${{\mathbb{F}}_{q}}$ then ord$\left( f \right)$ is the minimal period of this sequence.
If $f$ and $g$ are relative prime then ord$\left( fg \right)$ is the least common multiple of ord$\left( f \right)$ and ord$\left( g \right).$ Moreover, if ${{x}^{\ell }}-1$ is divisible by $f$ then it is also divisible by $\sigma \left( f \right).$ Therefore, the order of $f$ is divisible by the order of $\sigma \left( f \right)$ and so on. Therefore, $\text{ord}\left( f \right)=\text{ord}\left( \sigma \left( f \right) \right).$

\bigskip
\par\noindent{\bf  Theorem 5. } {\it  The minimal period of a periodic sequence over ${{\mathbb{F}}_{q}}$ is ${{q}^{N}}-1$ iff  its minimal polynomial is a product of distinct irreducible factors such that the degree of each factor divides $N$ and the least common multiple of orders of factors is exactly ${{q}^{N}}-1$.
}

 \bigskip
\par\noindent{\sl
Proof: } First we assume that $m={{m}_{1}}{{m}_{2}}\cdots {{m}_{h}}$ is the minimal polynomial of a periodic sequence such that ${{m}_{1}},{{m}_{2}},\cdots ,{{m}_{h}}$ are distinct irreducible polynomials of degrees dividing $N$. Then it follows from Theorem 2 that ${{q}^{N}}-1$ is a period of this sequence. Moreover, the minimal period of any periodic sequence is the order of its minimal polynomial. Here, the order of $m={{m}_{1}}{{m}_{2}}\cdots {{m}_{h}}$ is the least common multiple of orders of ${{m}_{1}},{{m}_{2}},\cdots ,{{m}_{h}}$  which is exactly ${{q}^{N}}-1$. 
Hence, the minimal period  is ${{q}^{N}}-1$.
Now assume that the minimal period of a periodic sequence over ${{\mathbb{F}}_{q}}$ is ${{q}^{N}}-1$. Then the minimal polynomial of this sequence divides ${{x}^{{{q}^{N}}-1}}-1.$ But ${{x}^{{{q}^{N}}}}-x$ itself is a square free product of irreducible polynomials of degrees dividing $N$ so the minimal polynomial is a product of distinct irreducible factors and the degree of each factor divides $N$. Moreover, the order of the minimal polynomial is ${{q}^{N}}-1$ and equal to the least common multiple of orders of its factors. The proof is now complete.
 
 \bigskip
\par\noindent{\bf Remark. } Theorem 5 is presented in \cite {Kyureghyan} in the case where $q=2$ without rigour proof.

 \bigskip
\par\noindent{\bf 	
Definition.}  An irreducible polynomial $f\in {{\mathbb{F}}_{q}}\left[ x \right]$ 
is called {\bf primitive } 
if $\text{ord}\left( f \right)={{q}^{\text{deg}\left( f \right)}}-1={{p}^{n\text{deg}\left( f \right)}}-1.$ Any root of a primitive polynomial of degree $N$ generates the cyclic group ${{\mathbb{F}}^*_{q^N}}$ so the number of primitive polynomials of degree $N$
is equal to the number of positive integers $\ell $ which are co-prime to $q$ and ${{q}^{N}}-1.$ For example, there are only 2 primitive polynomials of degree 3 over $\mathbb{F}_{2}$. They are ${{x}^{3}}+{{x}^{2}}+1$ and ${{x}^{3}}+x+1$. 
Since $\text{ord}\left( f \right)=\text{ord}\left( \sigma \left( f \right) \right)$ and $f$ is irreducible if and only if $\sigma \left( f \right)$ is irreducible, we have $f$ is primitive if and only if $\sigma \left( f \right)$ is primitive.
Recall that $k\left( f \right)$ denotes the smallest positive integer $\ell $ such that ${{\sigma }^{\ell }}\left( f \right)=f.$ We have $k\left( f \right)|n$ and $f$ is in ${{\mathbb{F}}_{{{p}^{k\left( f \right)}}}}\left[ x \right]$ so if $f$ is irreducible then $f$ divides ${{x}^{{{p}^{k\left( f \right)\text{deg}\left( f \right)}}-1}}-1.$ Therefore, if even more $f$ is primitive then it follows that $k\left( f \right)=n.$ 
Moreover, any root  of a  primitive polynomial $f\in {{\mathbb{F}}_{q}}[x]$ of degree $N$ is a generator element of the multiplicative group $\mathbb{F}_{{{q}^{N}}}^{*}$ so for  any index $i\in\{1,\cdots,q^N-2\}$ there is an another index $h\in\{1,\cdots,q^N-2\}$ such that $x^h-x^i-1$ is divisible by $f$.

\bigskip
\par\noindent{\bf   Theorem 6.} {\it 
If $f\in {{\mathbb{F}}_{p^n}}\left[ x \right]$ is primitive then  $g=f\sigma \left( f \right){{\sigma }^{2}}\left( f \right)\cdots {{\sigma }^{n-1}}\left( f \right)$ is in ${{\mathbb{F}}_{p}}\left[ x \right]$ and primitive over the prime field ${{\mathbb{F}}_{p}}.$
}

\bigskip
\par\noindent{\sl
Proof: }  We have seen that $g\in {{\mathbb{F}}_{p}}\left[ x \right]$ and it is irreducible over there. The order of $f$ is ${{p^n}^{\text{deg}\left( f \right)}}-1={{p}^{\text{deg}\left( g \right)}}-1$ which is exactly the order of $g$ over ${{\mathbb{F}}_{p}}$ so $g$ is primitive over ${{\mathbb{F}}_{p}}$. The proof is now complete.

 \bigskip
\par\noindent{\bf 	
Definition. } 
If the minimal polynomial of a periodic sequence is primitive of degree $N$ then we call this sequence a {\bf max-sequence}  of order $N$. 
Clearly, the minimal period of a max-sequence of order $N$ over ${{\mathbb{F}}_{q}}$ is ${{q}^{N}}-1.$
 It is well known that if $S=\left\{ {{s}_{0}},{{s}_{1}},\cdots  \right\}$ is a max-sequence of order $N$ then there is $\beta \in {{\mathbb{F}}_{{{q}^{N}}}^*}$ such that ${{s}_{j}}=\text{Tr}\left( {{\alpha }^{j}}\beta  \right)={{\alpha }^{j}}\beta +{{\left( {{\alpha }^{j}}\beta  \right)}^{q}}+{{\left( {{\alpha }^{j}}\beta  \right)}^{{{q}^{2}}}}+\cdots +{{\left( {{\alpha }^{j}}\beta  \right)}^{{{q}^{N-1}}}}$  for every $j=0,1,2\cdots .$               
Here, $\alpha$ denotes a root of the minimal polynomial of the max-sequence.
Indeed, it follows at once from the fact that the $q$th power of a determinant over the field ${{\mathbb{F}}_{q}}$ is giving by the following formula
\[{{\left| \begin{matrix}
   {{a}_{11}} & {{a}_{12}} & \cdots  & {{a}_{1N}}  \\
   {{a}_{21}} & {{a}_{22}} & \cdots  & {{a}_{2N}}  \\
   \vdots  & \vdots  & \ddots  & \vdots   \\
   {{a}_{N1}} & {{a}_{N2}} & \cdots  & {{a}_{NN}}  \\
\end{matrix} \right|}^{q}}=\left| \begin{matrix}
   a_{11}^{q} & a_{12}^{q} & \cdots  & a_{1N}^{q}  \\
   a_{21}^{q} & a_{22}^{q} & \cdots  & a_{2N}^{q}  \\
   \vdots  & \vdots  & \ddots  & \vdots   \\
   a_{N1}^{q} & a_{N2}^{q} & \cdots  & a_{NN}^{q}  \\
\end{matrix} \right|.\]
It is proved in \cite{Blackburn} that for a max-sequence of order $N$ and an index $i\in\{1,\cdots,q^N-2\}$ there is an another index $h\in\{1,\cdots,q^N-2\}$ such that 
$ S +{{\mathbb{E}}^{i}}\left( S \right)={{\mathbb{E}}^{h}}\left( S \right)$.
In fact, it follows at once from the fact that 
${{s}_{j}}+{{s}_{j+i}}=\text{Tr}\left( \left( 1+{{\alpha }^{i}} \right){{\alpha }^{j}}\beta  \right)=\text{Tr}\left( {{\alpha }^{h}}{{\alpha }^{j}}\beta  \right)={{s}_{j+h}},$
(because $\alpha$ is a generator element of the multiplicative group $\mathbb{F}_{{{q}^{N}}}^{*}$). On the other hand, note that in the trace representation 
${{s}_{j}}=\text{Tr}\left( {{\alpha }^{j}}\beta  \right)={{\alpha }^{j}}\beta +{{\left( {{\alpha }^{j}}\beta  \right)}^{q}}+{{\left( {{\alpha }^{j}}\beta  \right)}^{{{q}^{2}}}}+\cdots +{{\left( {{\alpha }^{j}}\beta  \right)}^{{{q}^{N-1}}}}$,  we have  $\beta \in {{\mathbb{F}}_{{{q}^{N}}}^*}$  so $\beta$ should be a power of $\alpha$ and consequently,  ${{s}_{j}}=\text{Tr}\left( {{\alpha }^{j+i_0}}  \right)$ for every $j=0,1,2\cdots .$   Therefore, {\bf max-sequences of the same minimal polynomial are equivalent.}  
The number of primitive polynomials of degree $N$ will gives us the number of non-equivalent max-sequences of order $N$. This number is equal to the number of positive integers $\ell $ which are co-prime to $q$ and ${{q}^{N}}-1.$ For example, there are only 2 primitive polynomials of degree 3 over $\mathbb{F}_{2}$. They are ${{x}^{3}}+{{x}^{2}}+1$ and ${{x}^{3}}+x+1$. Therefore, there are only two max-sequences of order 3 over  $\mathbb{F}_{2}$, which are 1110100 and 1110010.
Finally, we note that the trace operator Tr is a linear functional,  which  maps ${{\mathbb{F}}_{{{q}^{N}}}}$ as a vector space (over ${{\mathbb{F}}_{q}}$) onto ${{\mathbb{F}}_{q}}$.

\bigskip
\bigskip

\sect{DE BRUIJN SEQUENCES AND THEIR MODIFICATION}

\bigskip
\par\noindent 
A de Bruijn sequence of order $N $th over ${{\mathbb{F}}_{q}}$ is a ${{q}^{N }}-$periodic sequence $B=\left\{ {{b}_{0}},{{b}_{1}},\cdots  \right\}$ such that 
\[\left\{ \left( {{b}_{i}},{{b}_{i+1}},\cdots ,{{b}_{i+N -1}} \right):\text{ }i=0,1,\cdots ,{{q}^{N }}-1 \right\}=\mathbb{F}_{q}^{N }.\]
 Clearly, the minimal polynomial of any de Bruijn sequence is a power of $x-1$. Moreover, 
the weight of a  de Bruijn sequence $B$ is zero so its complexity $c\left( B \right)$ should satisfy ${{p}^{nN -1}}<c\left( B \right)<{{p}^{nN }}$ where ${{p}^{n}}=q.$
Let 
\[B\left( x \right)={{b}_{0}}+{{b}_{1}}x+\cdots +{{b}_{{{q}^{N }}-1}}{{x}^{{{q}^{N }}-1}}.\]
Then the minimal polynomial of $B$ is 
\[\frac{{{x}^{{{q}^{N }}}}-1}{\text{gcd}\left( {{x}^{{{q}^{N }}}}-1,B\left( x \right) \right)}={{\left( x-1 \right)}^{c\left( B \right)}}\]
so 
\[B\left( x \right)={{\left( x-1 \right)}^{{{q}^{N }}-c\left( B \right)}}f\left( x \right)\]
 where $f$ is co-prime with $x-1$ and $\text{deg}\left( f \right)=c\left( B \right)-1.$
We are interested in a modification of de Bruijn sequences. It is obtained from a  de Bruijn sequence by removing a zero from the unique allzero $N -$tuple. The minimal period of a  modified sequence of order $N$  is ${{q}^{N }}-1$. Hence, by Theorem 5 the minimal polynomial of a modified sequence is a product of distinct irreducible factors such that the degree of each factor divides $N $ and the least common multiple of orders of factors is exactly ${{q}^{N }}-1$.
On the other hand, the minimal polynomial of the modified sequence is
\[\frac{{{x}^{{{q}^{N }}-1}}-1}{\text{gcd}\left( {{x}^{{{q}^{N }}-1}}-1,B\left( x \right) \right)}=\frac{{{x}^{{{q}^{N }}-1}}-1}{\left( x-1 \right)\text{gcd}\left( {{x}^{{{q}^{N }}-1}}-1,f\left( x \right) \right)}=:m\left( x \right).\]
It is interesting to study in more details the irreducible factors of the minimal polynomial of a modified de Bruijn sequence. 

\bigskip
\par\noindent{\bf   Theorem 7.} {\it Any max-sequence is a modified de Bruijn sequence of the same order. Therefore, there exists a  de Bruijn sequence of any order.
The number of non-equivalent de Bruijn sequences of order $N$ is at least the number of primitive polynomials of the same order.
Recall that a periodic sequence is called a max-sequence if its  minimal polynomial is primitive.
}

\bigskip
\par\noindent{\sl
Proof: }  Let ${{s}_{j}}=\text{Tr}\left( {{\alpha }^{j+{{i}_{0}}}} \right)$ be a max-sequence of order $N$. Let $f$ be the  minimal polynomial of the max-sequence.  Assume that  
$\text{Tr}\left( {{\alpha }^{\ell -1}} \right)\ne 0$ and  $\text{Tr}\left( {{\alpha }^{\ell }} \right)=\text{Tr}\left( {{\alpha }^{\ell +1}} \right)=\cdots =\text{Tr}\left( {{\alpha }^{\ell +N-1}} \right)=0$. Then $\text{Tr}\left( {{\alpha }^{\ell +N}} \right)=-f\left( 0 \right)\text{Tr}\left( {{\alpha }^{\ell -1}} \right)\ne 0$ so the max-sequence is really a modified de Bruijn sequence of the  order  $N$.

\bigskip
\par\noindent 
Now we define ${{B}_{N}}$ the $q-$ary {\bf de Bruijn directed graph of order}  $N$ as follows.
Let the vertices of ${{B}_{N}}$ be the  elements of $\mathbb{F}_{q}^{N}$. Moreover, a vertex $\left( {{a}_{1}},{{a}_{2}},\cdots ,{{a}_{N}} \right)$ is directed to another vertex $\left( {{b}_{1}},{{b}_{2}},\cdots ,{{b}_{N}} \right)$ if $\left( {{a}_{2}},\cdots ,{{a}_{N}} \right)=\left( {{b}_{1}},{{b}_{2}},\cdots ,{{b}_{N-1}} \right).$ A de Bruijn sequence of order $N$ is the directed Hamilton cycle in this graph, which is certainly exists! Moreover, the diameter of this graph is $N$. Two non-equivalent de Bruijn sequences will give two distinct Hamilton cycles. For example, there are only two non-equivalent de Bruijn sequences of order 3 over  $\mathbb{F}_{2}$ are 11101000 and 11100010. 
Hence, the number of Hamilton cycles in ${{B}_{N}}$ is at least the number of primitive polynomials of order $N$. 
The number of directed Hamilton cycles in ${{B}_{N}}$ is at least the number of positive integers $\ell$ which are co-prime to $q$ and ${{q}^{N}}-1.$ Now we define the {\bf difference operator }
\[\mathbb{D}:\mathbb{F}_{q}^{N}\to \mathbb{F}_{q}^{N-1}\] by setting 
$\mathbb{D}\left( {{a}_{1}},\cdots ,{{a}_{N}} \right)=\left( {{a}_{2}}-{{a}_{1}},\cdots ,{{a}_{N}}-{{a}_{N-1}} \right).$ Then $\mathbb{D}$ induces an graph homomorphism from ${{B}_{N}}$ onto ${{B}_{N-1}}.$ Therefore, any de Bruijn sequence of order $N$ induces a de Bruijn sequence of order $N-1.$ For example,  two sequences 11101000 and 11100010 of order 3 induce the same sequence 0011 of order 2.

 \bigskip
\bigskip 

\sect{CONVOLUTION OF PERIODIC SEQUENCES}

\bigskip

\par\noindent  For two sequences $S=\left\{ {{s}_{0}},{{s}_{1}},\cdots  \right\}$ and $T=\left\{ {{t}_{0}},{{t}_{1}},\cdots  \right\}$ we define their convolution $S*T$ to be the sequence $U=\left\{ {{u}_{0}},{{u}_{1}},\cdots  \right\}$ by letting
\[{{u}_{\ell}}={{s}_{\ell}}{{t}_{0}}+{{s}_{\ell-1}}{{t}_{1}}+\cdots +{{s}_{0}}{{t}_{\ell}}\quad \text{   for  } \ell=0,1,2,\cdots .\]
For example, the convolution of $\left\{ 1,3,{{3}^{2}},\cdots  \right\}$ with itself over $\mathbb F_5$ is the Fibonacci sequence $\{1,1,2,3,\cdots\}$. 
Recall that 
$\mathcal M\left( f \right)$ denotes the vector space of periodic sequences having $f$ as a characteristic polynomial. Here $f\in {{\mathbb{F}}_{q}}\left[ x \right]$ is a monic polynomial.
The dimension of $\mathcal M\left( f \right)$ over $ {{\mathbb{F}}_{q}}$  is deg$\left( f \right).$	It is proved in \cite{Haukkanen}  that 
\[\mathcal{M}\left( f \right)*\mathcal{M}\left( g \right)\subseteq \mathcal{M}\left( fg \right)
\quad 
\text{ for all monic polynomials } f,g\in {{\mathbb{F}}_{q}}\left[ x \right].\]
This result is the best. For example, over $\mathbb F_5$ we have
\[\mathcal{M}\left( x-3 \right)*\mathcal{M}\left( x-3 \right)\ne \mathcal{M}\left( {{\left( x-3 \right)}^{2}} \right).\]
 Moreover, it is proved in \cite{Haukkanen2}  that
 \[\mathcal{M}\left( xf (x)\right)*\mathcal{M}\left( g \right)= \mathcal{M}\left( fg \right)
\quad 
\text{ for all monic polynomials } f,g\in {{\mathbb{F}}_{q}}\left[ x \right].\]
A typical example of convolution is Catalan number
\[{{c}_{n}}=\frac{1}{n+1}\left( \begin{matrix}
   2n  \\
   n  \\
\end{matrix} \right). \]
It is well known that 	\[{{c}_{\ell +1}}={{c}_{\ell }}{{c}_{0}}+{{c}_{\ell -1}}{{c}_{1}}+\cdots +{{c}_{0}}{{c}_{\ell }}\quad \text{   for  }\ell =0,1,2,\cdots \]
so in our notation $\mathbb{E}\left( C \right)=C*C$ where $C=\left\{ {{c}_{0}},{{c}_{1}},\cdots  \right\}.$

\bigskip
\bigskip
\par\noindent {\bf Acknowledgement. }
 Deepest appreciation is extended towards the NAFOSTED  (the National Foundation for Science and Techology Development in Vietnam) for the financial support.

\end{document}